\documentclass{article}

\begin{document}

\newtheorem{thm}{Theorem}[section]
\newtheorem{lem}[thm]{Lemma}
\newtheorem{rem}[thm]{Remark}
\newtheorem{cor}[thm]{Corollary}
\newtheorem{prop}[thm]{Proposition}
\def\proof{\noindent{\it Proof.\ }}
\newcommand{\bbfR}{{\rm I\kern-.1567em R}}
\newcommand{\R}{{\rm I\kern-.1567em R}}
\newcommand{\N}{{\rm I\kern-.1567em N}}
\newcommand{\Z}{{\sf Z\kern-.3567em Z}}
\newcommand{\ve}{\varepsilon}
\def\qed{\hfill$\Box$}
\newcommand{\be}{\begin{equation}} \newcommand{\ee}{\end{equation}}
\newcommand{\bea}{\begin{eqnarray}} \newcommand{\eea}{\end{eqnarray}}
\newcommand{\bean}{\begin{eqnarray*}} \newcommand{\eean}{\end{eqnarray*}}
\newcommand{\rf}[1]{(\ref {#1})}
\newcommand{\un}{\rm 1\!\!I}
\newcommand{\la}{\langle}
\newcommand{\ra}{\rangle}
\newcommand{\iom}[1]{\int_\Omega #1\, dx}
\newcommand{\idom}[1]{\int_{\partial\Omega} #1\; d\sigma}
\def\df{\stackrel{\rm df}=}
\def\t{\vartheta}
\def\div{\nabla\cdot}
\def\f{\varphi}
\def\r{\varrho}
\def\j{{\bf j}}
\def\h{{H^1(\Omega)}}
\font\eu=eufm10


\baselineskip=24pt



\title{Multiple solutions for equations involving bilinear, coercive and compact forms with 
applications to differential equations} 
\author{Robert Sta\'nczy\\
\small Instytut Matematyczny, Uniwersytet Wroc{\l}awski, \\
\small pl. Grunwaldzki 2/4, 50--384 Wroc{\l}aw, Poland.\\
\small {\tt stanczr@math.uni.wroc.pl}}

\date{\today}

\maketitle

\begin{abstract} 
The existence of multiple fixed points for the coercive, bilinear, 
compact forms defined in the cone in the Banach space. Multiple 
applications to the integral equations derived from BVPs for
differential euations are provided.
\end{abstract}




\section{Intoduction and motivation} 

It is well known that the quadratic equation
\be\label{qur}
u=au^2+u_0
\ee
can have either none, one or two solutions $u\in \R$, depending on 
the data $a>0$ and $u_0\in \R$. For example, if we assume that
\begin{equation}\label{qco}
4au_0<1
\end{equation}
than the existence of two positive solutions of \rf{qur} is guaranteed. 

In this paper we would like to show that this simple observation 
can be generalized if we replace quadratic term $ax^2$ 
with a bilinear form under suitable conditions. More 
specifically, we shall consider the equation in the cone $P$ in the 
Banach space $U$ with the norm $|\cdot|$ in the form
\be\label{ebf}
u=b(u,u)+u_0
\ee
for some given element $u_0\in P$ and bilinear, coercive and compact  
form $b$ defined on the product space $P\times P$. 

The assumption \rf{qco} guaranteeing the existence of two solutions 
for the quadratic equation \rf{qur} has to be adequately rephrased for \rf{ebf}  
as
\be\label{cbf}
4|b||u_0|<1\,
\ee
where $|b|$ denotes the norm of the bilinear norm.
However, the proof of this result we shall postpone to the next section.

Let us, however notice that the Banach Fixed Point Theorem for local 
contraction was used extensively to prove the existence of at least 
one fixed point (e.g. for the Navier-Stokes equations in \cite{CM} or for a Boltzmann equation) 
for the bilinear model equation like \rf{ebf}. This fixed 
point was located in the neighbourhood of $0$ thus making contraction approach feasible.

To prove the existence of two solutions we shall use 
another Fixed Point Theorem due to Krasnosielski (cf. \cite{Guo}) which allows to obtain 
more solutions if the nonlinear opearator has the required property
of ``crossing" identity twice, i.e. by the cone compression and the 
expansion on some appropriate convex subsets of the cone.

A quadratic nonlinearity being the simplest possible in the nonlinear world poses many questions about the global existence and uniqueness.
There are numerous models which posses such a structure, like the Navier Stokes euqation, the Boltzmann equation, the quadratic reaction diffusion equation (cf. 
\cite{HW}), 
the Smoluchowski coagulation equation or the system 
modelling chemotaxis and to name but a few. The problem of uniqueness of solutions for these equations attracted a lot of attention and only some
partial results are known. In some cases nonuniqueness occurs and the existence of two solutions can be proved. Sometimes one of the solution is a trivial one and
then the proof relies on finding a nontrivial one. In these dificult but important models one encounters another problem making our approach not feasible i.e. very 
common lack of compactness, thus if we would like to make our apprach feasible we are forced to consider some truncated baby model compatible with compact bilinear 
forms.

\section{The main result}

To prove the result announced in the previous section we shall 
use the following theorem \cite[Theorem 2.3.4]{Guo} originating from the works of Krasnosielski, cf., e.g., \cite{Kra}.
\begin{thm}\label{con}
Let $E$ be a Banach space, and let $P\subset E$ be a cone in $E$. Let $\Omega_1$ and
$\Omega_2$ be two bounded open sets in $E$ such that $0\in \Omega_1$ and $\overline{\Omega}_{1}\subset\Omega_2$. Let 
operator $A$ satisfy conditions
$$
\|Au\|\leq\|u\| {\rm \;\; for \;\; any \;\;} u\in P\cap\partial\Omega_1 {\rm \;\; and \;\;}
\|Au\|\geq\|u\| {\rm \;\; for \;\; any \;\;} u\in P\cap\partial\Omega_2 {\rm \;\;}
$$
or
$$
\|Au\|\geq\|u\| {\rm \;\; for \;\; any \;\;} u\in P\cap\partial\Omega_1 {\rm \;\; and \;\;}
\|Au\|\leq\|u\| {\rm \;\; for \;\; any \;\;} u\in P\cap\partial\Omega_2 {\rm \;\;}
$$
is satisfied. Then $A$ has at least one fixed point in $P\cap(\overline{\Omega}_2\setminus
\Omega_1)$.
\end{thm}

\begin{thm}
Assume that, for the given cone $P\subset E$, the bilinear and compact
form $b:P\times P\rightarrow P$ satisfies the following coercivity 
condition
\begin{equation}\label{coe}
\inf_{|u|=1,u\in P} b(u,u)>0.
\end{equation}
Then for any $u_0\in P$ as small as to satisfy \rf{cbf} the equation \rf{ebf} 
admits at least two solutions in P.
\end{thm}
\noindent
{\bf Proof.} Let us define the operator 
\begin{equation}\label{Tde}
Tu=b(u,u)+u_0
\end{equation}
then we shall apply Krasnosielski Theorem once as a cone-compression in the neighbouhodd of zero and secondly as a cone-expansion at infinity.

Notice that we have the following estimates
\begin{eqnarray}
\begin{array}{ll}
|Tu|\le |u_0|+B|u|^2,\\
|Tu|\ge |u_0|-B|u|^2,
\end{array}
\end{eqnarray}
where constant $B=|b|>0$ denotes the norm of the bilinear form $b$, i.e., the least (smallest??) constant $B$ satifying for any $u,v\in $ the inequality
$$
|b(u,v)|\le B|u||v|\,.
$$ 
for any $u,v\in E$.

Hence for sufficiently small $\rho_1>0$, i.e., such that $B\rho_1^2+\rho_1<|u_0|$ and any $u\in P$ and 
$|u|=\rho_1$ one has
\begin{equation}\label{1co}
|Tu|\ge |u_0|-B\rho_1^2 > \rho_1=|u|.
\end{equation}

Moreover, if we assume that there exists $\rho_2$ such that
\begin{equation}
|u_0|+B\rho_2^2<\rho_2
\end{equation}
then apparently for any $u\in P$ and $|u|=\rho_2$ one has
\begin{equation}\label{2co}
|Tu|\le |u_0|+B|u|^2<\rho_2=|u|.
\end{equation}
But this can be accomplished if we assume $B\rho_1^2+\rho_1< |u_0| < \rho_2-B\rho_2^2.$

Finally, for sufficiently large values of $\rho_3>0$ and any 
$u\in P$ and $|u|=\rho_3$, due to the coercivity assumption \rf{coe}, one has
\begin{equation}
b(u,u)\ge C|u|^2
\end{equation}
implying
\begin{equation}\label{3co}
|Tu|\ge C\rho_3^2 -|u_0|> \rho_3=|u|.
\end{equation}
To be more specific $\rho_3$ has to be so large that $|u_0|<C\rho_3^-\rho_3.$ 

Combining \rf{1co} with \rf{2co} we get that part of the cone P between 
the values $\rho_1$ and $\rho_2$ (in the $|\cdot|$ norm) is compressed 
while between the values $\rho_2$ and $\rho_3$ is expanded yielding the 
desired two fixed points in each set. Note that it might be necessary to 
distinguish between $\rho_2$ used in both sets as to prevent both fixed 
points to glue together. 

The last  the only assumption which should be made is to guarantee 
\rf{2co} to hold which follows readily from \rf{cbf}. 

\section{Examples of applications to differential equations}

\noindent
{\bf Example 1.}

Consider the following BVP for ODE for continuous postitive function $f$
\begin{eqnarray}
-u''(t)=(u(t))^2+f(t),\\
u(0)=u(1)=0.
\end{eqnarray}
This problem can be formulated in a required form
\begin{equation}
u=b(u,u)+u_0
\end{equation}
where the function
\begin{equation}
u_0(t)=\int_0^1 G(t,s)f(s)\,ds
\end{equation}
for the symetric, Green function given by $G(t,s)=t(1-s), 0\le t\le s\le 1$, while the bilinear form $b$ is defined, for any $u,v\in P$, by
$$b(u,v)=\int_0^1 G(t,s) u(s) v(s),$$
and the cone $P$ on which coercivity of $B$ holds can be defined as
$$P=\left\{u\ge 0:\inf_{t\in [a,b]} u(t) \ge \min\{a,1-b \} 
|u|_{\infty}\right\}\,.$$
For more applications to this kind of problems see \cite{PS, Lee}.

\noindent
Example 2.
Consider the following BVP for PDE, where $\Omega$ is an annulus in 
$\R^n$,
\begin{eqnarray}
\Delta u(x)=(u(x))^2+f(x)\,, \;\; x\in\Omega,\\
u|\partial \Omega =0.
\end{eqnarray}
Then the problem can be formulated as
\begin{equation}
u=b(u,u)+u_0
\end{equation}
where
\begin{equation}
u_0(x)=\int_\Omega G(x,y)f(y)\,ds
\end{equation}
for the appropriate, symmetric Green function $G$. The bilinear 
form $B$
is defined by
$B(u,v)=\int_0^1 G(x,y) u(y) v(y)$
and the cone $P$ on which coercivity of $B$ holds, can be expressed as $P\left\{u\ge 0:\inf_{x\in D} u(x) \ge \gamma |u|_{\infty}) \right\},$
for $\overline {D}\subset \Omega$ and some positive $\gamma$.
For more applications to this kind of problems see \cite{Ha}.

\noindent
{\bf Example 3.}
Consider for $u=u(x,t)$ the boundary value problem
\begin{equation}
u_t =\Delta u + u^2,\;\; u(x,0)=f(x)\,.
\end{equation}
Then one can formulate this problem as
\begin{equation}
u=b(u,u)+u_0
\end{equation}
where
\begin{equation}
u_0=S(t)f,
\end{equation}
while the bilinear form $b$
is defined by
\begin{equation}
b(u,v)=\int_0^t S(t-s)u(s)v(s)\,ds\,.
\end{equation}
Note that for any $t,s$ and the heat semigroup $S(t)$ one has
\begin{equation}
C|S(t-s)u(s)^2|_2\ge |S(t-s)u(s)^2|_1 \ge |u(s)^2|_1 \ge |u(s)|^2_2\,,
\end{equation}
hence after integration on $(0,t)$ and taking sup norm with respect to $t\in [0,T]$, the right hand side can be estimated by
\begin{equation}
\int_0^T |u(s)|^2_2 \,ds\,,
\end{equation}
while taking squared integral with respect to $t\in [0,T]$ the right hand side can be estimated by
\begin{equation}
T \left(\int_0^T |u(s)|^2_2 \right) \,ds\,.
\end{equation}
Moreover
\begin{equation}
|S(t)f|_2\le C|f|_2\,,
\end{equation}
hence the second condition for sup norm in $t$ follows. However one cannot guarantee the condition for coercivity. 
Note that due to \cite{HW} there are some results guaranteeing non uniqueness results. Our approach cannot yield this second solution, though
some nonuniqueness results are known for semilinear parabolic equations \cite{Wei}.

\noindent
{\bf Example 4. General case.}

Since the crucial assumptions in both examples which guarantee that the cone 
$P$ is invariant is the
following property of the Green function (or in general the kernel of Hammerstein operator involving the bilinar form $b(u,v)$) $G$
\begin{equation}\label{Ges}
\inf_{x\in U} G(x,y) \ge \gamma \sup_{x\in V}G(x,y)
\end{equation}
where $\gamma>0$ is independent of a set $V$ and its subset $U$ such that $\overline{U}\subset V$. Note that this 
condition holds either if $V$ is interval or annulus but we were not able to prove it for arbitrary domain, e.g. for a ball in higher dimension. If \rf{Ges} holds without any problems one can prove that properly defined cone is invariant under the action of the bilinear form, thus making it coercive. This problem can be illustrated in the second example, since if we replace an annulus with a ball the Green function is not bounded from above and also in the third one for the heat semigroup where the norms cannot be compared. This is the main obstacle in extending this kind of results to physically interesting models like the Navier Stokes equation, the Boltzmann equation or 
the Smoluchowski coagulation models sharing the property of quadratic nonlinearities.

\medskip\noindent\emph{Acknowledgments.\/} 
{\footnotesize This work has been partially supported by the Polish Ministry of Science project N N201 418839.}

\bigskip


\end{document}